\documentclass[a4paper]{article}
\usepackage[french]{babel}
\usepackage[latin1]{inputenc}
\usepackage[T1]{fontenc}
\usepackage{lmodern}
\usepackage{amsmath}
\usepackage{amsfonts}
\usepackage{stmaryrd}
\usepackage{graphicx}
\usepackage{graphics}
\usepackage{latexsym}
\usepackage{mathrsfs}
\usepackage{epsfig}
\usepackage{pgf}
\usepackage{tikz,tkz-tab}
\usepackage[all]{xy}
\usepackage{hyperref}
\tikzset { domaine/.style 2 args={domain=#1:#2} }
\usepackage{ncccomma}
\usepackage{changebar}
\usepackage{cancel}
\usepackage{manfnt}

\usetikzlibrary{snakes}
\usetikzlibrary{decorations}
\usetikzlibrary{calc,intersections,decorations.markings}

\newtheorem{defi}{Définition}[section]

\newtheorem{theo}[defi]{Th\'eor\`eme}
\newtheorem{lem}[defi]{Lemme}
\newcommand{\N}{\mathbb{N}}

\newcommand{\Q}{\mathbb{Q}}
\newcommand{\R}{\mathbb{R}}

\newcommand{\ud}{\mathrm{d}}
\newcommand{\Pre}{\mathfrak{p}}
\newcommand{\re}{\mathrm{Re}}
\usepackage{manfnt}

\makeatletter
\@addtoreset{footnote}{subsection}
\makeatother
\title{Théorème de Chebotarev effectif}
\author{Bruno Winckler\,\footnote{Université Bordeaux 1, Talence, bruno.winckler@math.u-bordeaux1.fr}}
\begin{document}

\maketitle

\noindent ABSTRACT: Let $K$ be a number field, and $L$ be a finite normal extension of $K$ with Galois group $G$. It is known that the number of Frobenius automorphisms corresponding to prime ideals, whose norms are less than $x$, is equivalent to the logarithmic integral as $x$ tends to infinity, and these automorphisms are well-distributed among the conjugacy classes of $G$: this is the Chebotarev theorem. The purpose of this paper is to compute the absolute constants of the error term appearing in a previous work about this theorem, due to Lagarias and Odlyzko.

\tableofcontents

\section{Introduction}

La répartition des nombres et idéaux premiers est un problème récurrent en théorie des nombres. Par exemple, on sait depuis le XIX${}^{e}$ siècle que la quantité de nombres premiers inférieurs à $x$ est asymptotiquement équivalente à $\frac{x}{\ln(x)}$, tandis que le théorème de la progression arithmétique de Dirichlet affirme que les nombres premiers se répartissent «~uniformément~» dans chaque classe de congruence possible. Le théorème de densité de Chebotarev, prouvé dans \cite{Che} en 1926, généralise les résultats de cette nature à l'ensemble des idéaux premiers d'une extension galoisienne $L/K$ de corps de nombres: les automorphismes de Frobenius associés aux idéaux premiers par la théorie algébrique des nombres se répartissent «~uniformément~» dans chaque classe de conjugaison du groupe de Galois. Ce théorème a de nombreuses applications, par exemple dans la théorie des courbes elliptiques ou des formes modulaires, comme l'illustre \cite{Ser}, et ceci suffit à justifier la recherche d'un développement asymptotique aussi précis que possible. Dans cette direction, \cite{Lag} donne les deux premiers termes du développement asymptotique de la fonction de décompte de ces idéaux premiers à des constantes non explicites près, et \cite{Oes} annonce une explicitation de ces constantes si on admet l'hypothèse de Riemann généralisée. En vue d'applications diophantiennes qui feront l'objet d'une étude ultérieure, on cherche à exhiber de telles constantes, et c'est l'objet de cet article. On précise le contexte ci-dessous.\\

Soit $L/K$ une extension galoisienne de corps de nombres, $d_L$ le discriminant absolu de $L$, $n_L$ le degré de $L$ sur $\Q$. Si $G$ est le groupe de Galois de cette extension, on désigne un ensemble de $G$ stable par conjugaison par la lettre $C$; pour tout $x > 1$, la fonction $\pi_C(x)$ décompte le nombre d'idéaux premiers $\mathfrak{p}$ de $K$ de norme inférieure ou égale à $x$, qui ne se ramifient pas dans $L$, et tels que $\left[\frac{L/K}{\Pre}\right] \in C$, où $\left[\frac{L/K}{\Pre}\right]$ est l'unique automorphisme de Galois (à conjugaison près) qui, réduit modulo $\mathfrak{p}$, coïncide avec l'automorphisme de Frobenius.

Si $C = G$, alors $\pi_C(x) = \pi_K(x) \sim \frac{x}{\ln(x)}$ quand $x\to \infty$. Le théorème de Chebotarev énonce que plus généralement, on a \[\pi_C(x) \sim \frac{|C|}{|G|}\frac{x}{\ln(x)}\textrm{ quand }x \to \infty.\]Plus précisément, si $\mathrm{Li}(x) = \int_2^x \frac{\ud t}{\ln(t)}$ et $\zeta_L(s) = \prod\limits_{\Pre \in \mathrm{Spec}(\mathrm{O}_L)\setminus\{0\}} \left(1-\mathrm{N}(\Pre)^{-s}\right)^{-1}$, alors on a une forme effective, qui s'exprime à l'aide de l'éventuel zéro positif $\beta$ de la fonction $\zeta_L$ tel que $0<1-\beta \leq \frac{1}{4\ln(d_L)}$.
\begin{theo}[Théorème de Chebotarev effectif]\label{chebotarev}Conservons les notations ci-dessus. On a, pour tout $x \geq \exp\left(8n_L (\ln(150867d_L^{44/5}))^2\right)$,\[\left|\pi_C(x) - \frac{|C|}{|G|}\mathrm{Li}(x)\right| \leq \frac{|C|}{|G|}\mathrm{Li}(x^\beta)+C_0x \exp\left(-\frac{1}{99} \sqrt{\frac{\ln(x)}{n_L}}\right),\]avec $C_0 = 783846699796966 < 7,84\cdot 10^{14}$.\end{theo}
On en sait un peu plus sur $\beta$: d'après \cite{Stark}, page 148, il existe une constante absolue et effectivement calculable $c$ telle que $\beta < 1-\frac{1}{cd_L^{1/n_L}}$; d'après \cite{BG}, il est fort probable que $c = \frac{\pi^2}{6}$ convienne. Le terme d'erreur est nettement plus petit si on suppose que l'hypothèse de Riemann est vérifiée pour $\zeta_L$:

\begin{theo}\label{chebotarevGRH}Supposons que l'hypothèse de Riemann généralisée est vraie pour $\zeta_L$. Alors, pour tout $x\geq 2$,
\begin{eqnarray}\left|\pi_C(x) - \frac{|C|}{|G|}\mathrm{Li}(x)\right| & \leq & \frac{|C|}{|G|}\sqrt{x}\Big[\left(32+\frac{181}{\ln(x)}\right)\ln(d_L)\nonumber\\
& \quad + &\left(28\ln(x)+330+\frac{1655}{\ln(x)}\right)n_L\Big]. \nonumber\end{eqnarray}\end{theo}

\emph{Remarque.} Dans \cite{Oes}, Oesterlé énonce même, sans démonstration,\[\left|\pi_C(x) - \frac{|C|}{|G|}\mathrm{Li}(x)\right| \leq \frac{|C|}{|G|}\sqrt{x}\left(\ln(d_L)\left(\frac{1}{\pi}+\frac{5,3}{\ln(x)}\right)+n_L\left(2+\frac{\ln(x)}{2\pi}\right)\right). \]
Ces deux théorèmes sont démontrés dans \cite{Lag}, mais avec des constantes non explicites. Pour trouver ces constantes, on suit la stratégie de démonstration proposée par \cite{Lag}, en soignant les majorations.

On remercie F. Pazuki, P. Autissier pour leurs conseils avisés, et K. Belabas pour l'intérêt manifesté pour ces résultats.

\section{Fil de la preuve}

L'essentiel de la preuve consiste en la recherche d'une formule asymptotique avec un terme d'erreur explicite pour $\psi_C(x) = \psi_C(x,L/K)$, une sorte de fonction pondérée de décompte des puissances d'idéaux premiers, intimement liée à $\pi_C(x,L/K)$. Cette fonction classique est définie par
\[\psi_C(x,L/K) = \sum_{\substack{\mathrm{N}_{K/\Q}(\Pre)^m \leq x \\ \Pre \textrm{ non ramifié} \\ \left[\frac{L/K}{\Pre}\right]^m = C}} \ln(\mathrm{N}_{K/\Q}(\Pre)).\]

Conformément à \cite{Lag}, voici les étapes de la preuve:
\begin{enumerate}
\item La fonction $\psi_C$ est très proche d'une transformée de Mellin inverse tronquée: \[I_C(x,T) = \frac{1}{2\pi i}\int_{\sigma_0-iT}^{\sigma_0+iT} F_C(s) \frac{x^s}{s}\ud s,\]où la fonction $F_C$ est construite en faisant une combinaison linéaire de dérivées logarithmiques de fonctions $\mathrm{L}$, en profitant de la formule d'orthogonalité des caractères de $G$ pour exclure tous les idéaux premiers qui ne nous intéressent pas; c'est dans le même esprit que la preuve traditionnelle du théorème de la progression arithmétique de Dirichlet est construite (voir \cite{Ser2} par exemple). On cherche à étudier l'écart $R_1(x,T)$ entre $\psi_C$ et $I_C$. Cet écart est explicité dans la section 3.
\item La fonction $F_C$ peut, en fait, être écrite comme combinaison linéaire de dérivées logarithmiques de fonctions $\mathrm{L}$ de Hecke. Toutes les singularités de $F_C$, qui sont des pôles simples, apparaissent en tant que zéros et pôles de $\zeta_L$. Cette étape est l'objet de la section 4, et est déjà entièrement explicite dans \cite{Lag}: je citerai le résultat en question \emph{sans le redémontrer}.
\item L'intégrale $I_C(x,T)$ diffère d'une intégrale sur un contour $B_C(x,T)$ par un terme d'erreur $R_2(x,T)$ (la fameuse étape où on «~translate à gauche la droite ( ! ) d'intégration~»). On a ici besoin de résultats de densité sur les zéros non triviaux de $\zeta_L$ pour estimer $R_2$; c'est l'objet de la section 4.
\item L'intégrale sur un contour est évaluée grâce au théorème des résidus, dans la section 6. L'intégrande a pour pôles les zéros et les pôles de $\zeta_L$; on en tire un terme principal $\frac{|C|}{|G|}x$ (c'est le pôle $s=1$) et une somme $S(x,T)$ indicée par les zéros de $\zeta_L$ à l'intérieur du contour de $B_C(x,T)$. On en déduit une formule explicite tronquée pour $\psi_C(x)$ avec un terme d'erreur inconditionnel, qui est donnée dans le théorème \ref{formulexplicite} de la section 6.
\item Il faut à présent étudier la somme sur les zéros $S(x,T)$. Si on suppose vraie l'hypothèse de Riemann généralisée, on peut directement aller à la section 8. Sinon, on a besoin de mieux connaître l'emplacement des zéros non triviaux: on dévoile des régions proches de $\re(s)=1$ sans zéros pour $\zeta_L$, dans la section 7.
\item La formule asymptotique $\psi_C(x) \sim \frac{|C|}{|G|}x$ avec un reste explicité en découle, par un bon choix de $T$ en fonction de $x$, visant à minimiser les termes d'erreur. C'est traité dans la section 8.
\item Enfin, la formule explicite pour $\pi_C$ s'obtient par une transformée d'Abel, à la fin de ce papier.
\end{enumerate}

Les résultats non prouvés ici sont ceux qui ont déjà été explicités dans \cite{Lag}, auquel cas on renvoie le lecteur à ce papier. Dans un souci d'auto-suffisance, on rappelle néanmoins toutes les définitions nécessaires à la compréhension de la preuve. Enfin, quand on remplace les constantes réelles par des nombres rationnels, cela se fait systématiquement grâce aux fractions continuées réduites à l'ordre 3.

\section{Fonctions $\mathrm{L}$ d'Artin et transformées de Mellin}

Par commodité, on note $\pi_C(x) = \pi_C(x,L/K)$, $\psi_C(x) = \psi_C(x,L/K)$ et $\mathrm{N} = \mathrm{N}_{K/\Q}$. Enfin, les sommes $\sum\limits_\Pre$ sont toujours indicées par les idéaux premiers non nuls de l'anneau des entiers de $K$.

L'article \cite{Lag} prouve que la fonction $F_C$ citée dans la section précédente a pour expression \[F_C(s) = \sum_{\Pre} \sum_{m=1}^\infty \theta(\Pre^m) \ln(\mathrm{N}(\Pre))(\mathrm{N}(\Pre))^{-ms},\]où pour $\Pre$ non ramifié dans $L$ on a $\theta(\Pre^m) = 1$ si $\left[\frac{L/K}{\Pre}\right]^m = C$, $\theta(\Pre^m) = 0$ sinon, et $|\theta(\Pre^m)| \leq 1$ si $\Pre$ se ramifie dans $L$.

On voit alors que si on oublie un instant les facteurs correspondant aux idéaux premiers ramifiés, $\psi_C(x)$ est une somme partielle des coefficients de $F_C(s)$. Pour obtenir $\psi_C(x)$ à partir de $F_C(s)$, on utilise la version tronquée de la transformée inverse de Mellin:
\begin{lem}Si $y > 0$, $\sigma >0$ et $T>0$, alors 
\begin{eqnarray}
\left|\frac{1}{2\pi i}\int_{\sigma-iT}^{\sigma+iT} \frac{y^s}{s}\ud s - 1\right| & \leq & y^\sigma \min(1,T^{-1}|\ln(y)|^{-1}) \textrm{ si }y>1,\nonumber \\
\left|\frac{1}{2\pi i}\int_{\sigma-iT}^{\sigma+iT} \frac{y^s}{s}\ud s - \frac{1}{2}\right| & \leq & \sigma T^{-1} \textrm{ si }y=1,\nonumber \\
\left|\frac{1}{2\pi i}\int_{\sigma-iT}^{\sigma+iT} \frac{y^s}{s}\ud s\right| & \leq & y^\sigma \min(1,T^{-1}|\ln(y)|^{-1}) \textrm{ si }0<y<1.\nonumber 
\end{eqnarray}\end{lem}

\emph{Démonstration.} Voir \cite{Dav}, pages 109--110, par exemple.\,$\square$\\

Soient $\sigma_0>1$ et $x\geq 2$. On définit $I_C(x,T) = \frac{1}{2\pi i}\int_{\sigma_0-iT}^{\sigma_0+iT} F_C(s) \frac{x^s}{s}\ud s$. Comme la série de Dirichlet définissant $F_C$ est absolument convergente pour $\re(s)>1$, on peut intégrer terme à terme (suivant le lemme 3.1) pour obtenir
\begin{eqnarray}\label{psic}\left|I_C(x,T) - \sum_{\substack{\Pre, m\\ \mathrm{N}(\Pre)^m \leq x}} \theta(\Pre^m)\ln(\mathrm{N}(\Pre))\right| & \leq & \sum_{\substack{\Pre, m \\ \mathrm{N}(\Pre^m) = x}} (\ln(\mathrm{N}(\Pre)) + \sigma_0 T^{-1}) \nonumber\\ & + & R_0(x,T),\end{eqnarray}
où
\begin{equation}\label{defR0}R_0(x,T) = \sum_{\substack{\Pre, m \\ \mathrm{N}(\Pre^m) < x}} \left(\frac{x}{\mathrm{N}(\Pre^m)}\right)^{\sigma_0} \min\left(1, T^{-1}\left|\ln\left(\frac{x}{\mathrm{N}(\Pre^m)}\right)\right|^{-1}\right)\ln(\mathrm{N}(\Pre)).\end{equation}
La somme du membre de gauche dans (\ref{psic}) égale $\psi_C(x)$, aux termes ramifiés près. On a $\mathrm{N}(\Pre) \geq 2$, et d'après \cite{Ser} (prop. 5 du n°1.3) on a \[\sum_{\Pre \textrm{ ramifié}} \ln(\mathrm{N}(\Pre)) \leq \frac{2}{|G|}\ln(d_L). \]Donc,
\begin{eqnarray}
\left|\sum_{\substack{\Pre, m\\ \mathrm{N}(\Pre^m) \leq x}} \theta(\Pre^m) \ln(\mathrm{N}(\Pre)) - \psi_C(x) \right| & \leq & \sum_{\substack{\Pre, m\\ \Pre \textrm{ ramifié} \\ \mathrm{N}(\Pre^m) \leq x}} \ln(\mathrm{N}(\Pre)) \nonumber \\ & \leq & \sum_{\Pre\textrm{ ramifié}} \ln(\mathrm{N}(\Pre)) \sum_{\substack{m \\ \mathrm{N}(\Pre^m)\leq x}} 1 \nonumber \\
& \leq & \frac{\ln(x)}{\ln(2)}\sum_{\Pre \textrm{ ramifié}} \ln(\mathrm{N}(\Pre)) \nonumber \\ & \leq & \frac{2}{\ln(2)}\frac{\ln(x)\ln(d_L)}{|G|}. \label{psic2}
\end{eqnarray}
Comme il y a au plus $n_K$ paires distinctes $(\Pre,m)$ telles que $\mathrm{N}(\Pre^m) = x$, on a 
\begin{equation}\sum_{\substack{\Pre, m \\ \mathrm{N}(\Pre^m) = x}} \ln(\mathrm{N}(\Pre)) \leq n_K \ln(x) \textrm{ et } \sum_{\substack{\Pre, m \\ \mathrm{N}(\Pre^m) = x}} \sigma_0 T^{-1} \leq n_K \sigma_0 T^{-1}. \label{psic3}\end{equation}
Alors, (\ref{psic}), (\ref{psic2}) et (\ref{psic3}) entraînent \[\psi_C(x) = I_C(x,T) + R_1(x,T),\]
où
\begin{equation}R_1(x,T) \leq \frac{2}{\ln(2)}\frac{\ln(x)\ln(d_L)}{|G|} + n_K\sigma_0 T^{-1} + n_K \ln(x) + R_0(x,T).\end{equation}
Le reste de la section est voué à estimer $R_0(x,T)$. À présent, on prend $\sigma_0 = 1+\frac{1}{\ln(x)}$ (on a alors $x^{\sigma_0} =ex$). Écrivons $R_0(x,T) = S_1+S_2+S_3$, où $S_1$ est la somme dans (\ref{defR0}) n'impliquant que les idéaux $\Pre$ tels que $|\mathrm{N}(\Pre^m) -x|\geq \frac{1}{4}x$, tandis que $S_2$ est indicé par ceux tels que $|x - \mathrm{N}(\Pre^m)| \leq 1$, et $S_3$ récupère les termes restants.

Pour $\Pre$ tel que $|\mathrm{N}(\Pre^m) -x|\geq \frac{1}{4}x$, on a $\min\left(1,T^{-1}\left|\ln\left(\frac{x}{\mathrm{N}(\mathfrak{p})^m}\right)\right|^{-1}\right) \leq \frac{T^{-1}}{\ln\left(\frac{5}{4}\right)}$. On obtient alors,\[S_1 \leq \frac{xT^{-1}e}{\ln\left(\frac{5}{4}\right)}\sum_{\mathfrak{p},m} \mathrm{N}(\mathfrak{p})^{-m\sigma_0} \ln(\mathrm{N}(\Pre)) = \frac{xT^{-1}e}{\ln\left(\frac{5}{4}\right)}\left(-\frac{\zeta_K'}{\zeta_K}(\sigma_0)\right).\]
La dernière majoration dépend de la comparaison entre $-\frac{\zeta'_\Q}{\zeta_\Q}(\sigma_0)$ et $(\sigma_0-1)^{-1}$. De l'analyse classique permet d'obtenir:
\begin{lem}On a, pour $\sigma>1$, \[-\frac{\zeta'_K}{\zeta_K}(\sigma) \leq -n_K\frac{\zeta'_\Q}{\zeta_\Q}(\sigma) \leq n_K(\sigma-1)^{-1}.\]\end{lem}

\emph{Démonstration.} La première inégalité est prouvée dans \cite{Lag}, page 426, et la deuxième inégalité se déduit de l'égalité \begin{equation}\label{zeta}\zeta_\Q(\sigma) = \frac{\sigma}{\sigma-1}-\sigma I(\sigma),\end{equation}où $I(\sigma) =\int_1^\infty (t-[t])t^{-\sigma-1}\ud t$. Sachant cela, on obtient facilement\[\frac{\zeta_\Q'}{\zeta_\Q}(\sigma)+\frac{1}{\sigma-1} > \frac{1-(2\sigma-1)I(\sigma)}{\zeta(\sigma)(\sigma-1)}\]pour tout $\sigma> 1$. Comme $I(\sigma) = \frac{1}{\sigma}\left(\frac{\sigma}{\sigma-1}-\zeta_\Q(\sigma)\right)$, l'inégalité à démontrer se réduit à l'inégalité $\zeta_\Q(\sigma) \geq \frac{\sigma^2}{(\sigma-1)(2\sigma-1)}$, qui est un exercice simple d'analyse
.\,$\square$\\

Reprenons. On a\begin{equation}\label{S1}S_1 \leq \frac{xT^{-1}e}{\ln\left(\frac{5}{4}\right)}\left(-n_K\frac{\zeta_\Q'}{\zeta_\Q}(\sigma_0)\right) \leq \frac{e}{\ln\left(\frac{5}{4}\right)}n_KxT^{-1}\ln(x).\end{equation}
La majoration de $S_2$ est poussée plus loin que dans \cite{Lag} assez facilement: on a, pour $x\geq 2$,
\begin{equation}\label{S2}S_2 \leq 2n_K \ln(x+1)\left(\frac{x}{x-1}\right)^{\sigma_0} \leq \frac{4e\ln(3)}{\ln(2)}n_K\ln(x).\end{equation}
Enfin, pour $S_3$ (qui est indicé par les idéaux $\Pre$ tels que $1<|N(\Pre^m)-x|< \frac{x}{4}$) on a, en vertu de l'inégalité $\left|\ln\left(\frac{x}{n}\right)\right|^{-1} \leq \frac{2n}{|x-n|}$ (pour $2n \geq x$),\begin{eqnarray}
S_3 & \leq & \frac{4e^2}{3^{1+1/\ln(2)}}T^{-1}\ln\left(\frac{5x}{4}\right)\sum_{\substack{n \\ 1 < |n-x| < \frac{x}{4}}} \left|\ln\left(\frac{x}{n}\right)\right|^{-1}\sum_{\substack{\Pre,m \\ \mathrm{N}(\Pre)^m = n}} 1\nonumber\\
    & \leq & \frac{8e^2}{3^{1+1/\ln(2)}}n_KT^{-1}\ln\left(\frac{5x}{4}\right)\sum_{\substack{n \in \N \\ 1 < |n-x| < \frac{x}{4}}} \frac{n}{|x-n|}  \nonumber\\
    & \leq & \frac{8e^2}{3^{1+1/\ln(2)}}n_KT^{-1}\ln\left(\frac{5x}{4}\right)\sum_{\substack{k \in \N \\ 1 < k < \frac{x}{4}}} \left(1+\frac{x}{k}\right) \nonumber\\
    & \leq &  \frac{8e^2}{3^{1+1/\ln(2)}}n_KT^{-1}x\ln\left(\frac{5x}{4}\right)\left(\frac{1}{4}+\ln\left(\frac{x}{2}\right)\right) \nonumber \\
    & \leq & \frac{8e^2}{3^{1+1/\ln(2)}}n_KT^{-1}x(\ln(x))^2 \label{S3}
\end{eqnarray}
En regroupant (\ref{S1}), (\ref{S2}) et (\ref{S3}), on obtient une majoration de $R_0$, puis de $R_1$:
\begin{eqnarray}
R_0(x,T) &\leq & \frac{4e\ln(3)}{\ln(2)}n_K \ln(x) +\left(\frac{8e^2}{3^{1+1/\ln(2)}}+\frac{e}{\ln\left(\frac{5}{4}\right)\ln(2)}\right)n_K T^{-1} x (\ln(x))^2 \nonumber \\
         &\leq & \frac{69}{4}n_K \ln(x)+\frac{65}{3}n_K T^{-1} x (\ln(x))^2\end{eqnarray}
valable pour tout $x \geq 2$ et tout $T \geq 1$, et donc,
\begin{equation}\label{majR1}R_1(x,T) \leq \frac{2}{\ln(2)}\frac{\ln(x)\ln(d_L)}{|G|}+\frac{73}{4}n_K\ln(x) + \frac{145}{6}n_KT^{-1}x(\ln(x))^2.\end{equation}

\section{Réduction au cas des fonctions $\mathrm{L}$ de Hecke}

Soit $g \in G$, et soit $H$ le sous-groupe de $G$ engendré par $g$, $E$ le corps fixé par $H$, et notons avec la lettre $\chi$ les caractères irréductibles de $H$ (qui sont de dimension 1).
\begin{lem}On a \[F_C(s) = -\frac{|C|}{|G|}\sum_\chi \bar{\chi}(g) \frac{\mathrm{L}'}{\mathrm{L}}(s,\chi,L/E).\]\end{lem}

\emph{Démonstration.} Voir \cite{Lag}, pages 429--431.\,$\square$\\

Ainsi, $F_C$ s'écrit non seulement à l'aide de fonctions $\mathrm{L}$ d'Artin correspondant à des caractères non linéaires de $G$, mais on peut aussi les exprimer comme combinaison linéaire de dérivées logarithmes de fonctions $\mathrm{L}$ de Hecke \emph{abéliennes}. L'article \cite{Lag} utilise ceci pour améliorer la localisation et la densité des singularités de $F_C$. 


On a donc
\begin{equation}\label{rappelIC}I_C(x,T) = -\frac{|C|}{|G|} \sum_\chi \bar{\chi}(g)\frac{1}{2\pi i}\int_{\sigma_0-iT}^{\sigma_0+iT} \frac{x^s}{s}\frac{\mathrm{L}'}{\mathrm{L}}(s,\chi,L/E)\ud s,\end{equation}et on a vu que $\psi_C$ s'écrit comme somme de cette intégrale et de $R_1$, qu'on a déjà estimé en (\ref{majR1}). Il s'agit donc, à présent, d'évaluer les intégrales présentes dans (\ref{rappelIC}). D'où la nécessité de borner le nombre de singularités de $\frac{\mathrm{L}'}{\mathrm{L}}$. C'est l'occasion de rappeler quelques propriétés classiques de ces fonctions $\mathrm{L}$, qu'on peut retrouver dans le chapitre 5 de \cite{Neu}.

Comme les corps $L$ et $E$ sont fixés, on peut omettre $L/E$ dans l'écriture de $\mathrm{L}(s,\chi,L/E)$. Soit $F(\chi)$ le conducteur de $\chi$, et soit
\[A(\chi) = d_E \mathrm{N}_{E/\Q} (F(\chi))\label{defAchi}.\]
On définit $\delta(\chi)$ comme étant égal à 1 pour le caractère principal, et 0 sinon. Enfin, pour rappel, pour chaque $\chi$ il existe deux entiers naturels $a(\chi)$ et $b(\chi)$ tels que $a(\chi) + b(\chi) = n_E$, de sorte que si on pose \[\gamma_\chi(s) = \left(\pi^{-(s+1)/2} \Gamma\left(\frac{s+1}{2}\right)\right)^{b(\chi)} \left(\pi^{-s/2} \Gamma\left(\frac{s}{2}\right)\right)^{a(\chi)}\] 
et 
\begin{equation}\label{defxi}\xi(s,\chi) = (s(s-1))^{\delta(\chi)} A(\chi)^{s/2} \gamma_\chi(s) L(s,\chi),\end{equation}
alors on a l'équation fonctionnelle 
\[\xi(1-s,\bar{\chi}) = W(\chi)\xi(s,\chi),\]
où $W(\chi)$ est une constante de module 1. De plus, $\xi(\cdot,\chi)$ est une fonction entière d'ordre 1 qui ne s'annule pas en 0, donc 
\begin{equation}\label{hadamardproduct}\xi(s,\chi) = e^{B_1(\chi)+B(\chi)s} \prod_\rho \left(1-\frac{s}{\rho}\right) e^{s/\rho}\end{equation}
pour des constantes $B_1(\chi)$ et $B(\chi)$, où $\rho = \beta + i \gamma$ parcourt l'ensemble des zéros (non triviaux) de $L(s,\chi)$ tels que $0<\beta<1$. La lettre $\rho$ désignera toujours ces zéros. En dérivant logarithmiquement (\ref{defxi}) et (\ref{hadamardproduct}), on obtient la relation importante
\begin{eqnarray}\label{identiteL'L}\frac{\mathrm{L}'}{\mathrm{L}}(s,\chi) = B(\chi) + \sum_\rho \left(\frac{1}{s-\rho}+\frac{1}{\rho}\right) - \frac{1}{2}\ln(A(\chi)) - \delta(\chi)\left(\frac{1}{s}+\frac{1}{s-1}\right) - \frac{\gamma_\chi'}{\gamma_\chi}(s),\end{eqnarray}
valable pour tout nombre complexe $s$ où ces quantités sont définies. On ne sait pas déterminer exactement la dépendance de $B$ en fonction de $\chi$, mais on sait tout de même déduire de l'équation fonctionnelle ceci:
\begin{lem}En conservant les notations ci-dessus, \[\re(B(\chi)) = - \sum_\rho \re\left(\frac{1}{\rho}\right),\]et
\begin{eqnarray}\frac{\mathrm{L}'}{\mathrm{L}}(s,\chi)+\frac{\mathrm{L}'}{\mathrm{L}}(s,\bar{\chi}) &= & \sum_\rho \left(\frac{1}{s-\rho}+\frac{1}{s-\bar{\rho}}\right)-\ln(A(\chi)) - 2\delta(\chi)\left(\frac{1}{s}+\frac{1}{s-1}\right) \nonumber \\ & \quad - & 2\frac{\gamma_\chi'}{\gamma_\chi}(s).\label{identiteL'L2}\end{eqnarray}\end{lem}

\emph{Démonstration.} Voir \cite{Odl}.\,$\square$\\

On établit ici quelques résultats préliminaires.
\begin{lem}Si $\re(s) > 1$, alors $\left|\frac{\mathrm{L}'}{\mathrm{L}}(s,\chi)\right| \leq \frac{n_E}{\re(s)-1}$.\label{L'surL}\end{lem}

\emph{Démonstration.} En comparant leurs séries de Dirichlet, on voit que \[\left|\frac{\mathrm{L}'}{\mathrm{L}}(s,\chi) \right| \leq -\frac{\zeta'_E}{\zeta_E}(\re(s)),\]et on en arrive aisément au résultat du lemme.\,$\square$

\begin{lem}\begin{enumerate}
\item Si $\re(z)\geq a$, avec $a\geq 1$, alors \[\left|\frac{\Gamma'}{\Gamma}(z) \right|\leq \ln(|z|)+\frac{\pi}{2}+\frac{1}{a}.\]
\item Si $|\mathrm{Im}(z)| \geq b \geq 1$, alors \[\left|\frac{\Gamma'}{\Gamma}(z) \right|\leq \ln(|z|)+\pi\left(1+\frac{1}{2b}\right)+\frac{1}{2b}.\]
\item Si $|z+k| \geq \frac{1}{8}$ pour tout entier naturel $k$, alors\[\left|\frac{\Gamma'}{\Gamma}(z) \right|\leq \ln(|z|)+\frac{83}{5}.\]\end{enumerate}\label{ineqdigamma}\end{lem}

\emph{Démonstration.} Le produit de Weierstrass de $\Gamma$ nous enseigne, en considérant sa dérivée logarithmique, que \[-\frac{\Gamma'}{\Gamma}(z) = \gamma_0 + \sum_{n=1}^\infty \left(\frac{1}{z+n-1}-\frac{1}{n}\right)\]pour tout complexe $z$ différent des entiers négatifs ($\gamma_0$ désigne la constante d'Euler-Mascheroni). Alors, en utilisant la formule d'Euler-McLaurin (pour $z > 0$) avec cette somme, on obtient que
\[\label{reprintgamma}\frac{\Gamma'}{\Gamma}(z) = \ln(z)-\frac{1}{2z}+\int_0^\infty \frac{B_1(\{x\})}{(z+x)^2}\ud x,\]
où $B_1 = X-\frac{1}{2}$ est le premier polynôme de Bernoulli; l'égalité pour tout $z$ complexe (à l'exception de $z\leq 0$) s'obtient par le principe du prolongement analytique. On en déduit que, pour $\re(z)\geq a \geq 1$,
\[\left|\frac{\Gamma'}{\Gamma}(z)\right| \leq \ln(|z|)+\frac{\pi}{2}+\frac{1}{a}.\]
On procède de même pour l'inégalité dans le domaine $|\mathrm{Im}(z)|\geq 1$. Passons à la troisième inégalité: si $\re(z)<1$ et $|\mathrm{Im}(z)| < 1$ (notons que le cas $|\mathrm{Im}(z)|\geq 1$ est déjà traité et vérifie l'inégalité annoncée), l'équation fonctionnelle vérifiée par $\Gamma$ permet de montrer que 
\[\frac{\Gamma'}{\Gamma}(z) = \frac{\Gamma'}{\Gamma}(z+m) - \sum_{k=0}^{m-1} \frac{1}{z+k}\]
pour tout entier $m>0$. On pose $m = -\left[\re(z)-1\right]$, de sorte que $\re(z+m)\geq 1$ (c'est le plus petit entier naturel à vérifier cette inégalité), ce qui permet de borner $\frac{\Gamma'}{\Gamma}(z+m)$:\[\left|\frac{\Gamma'}{\Gamma}(z+m)\right| \leq \ln\left(\sqrt{5}\right)+ \frac{\pi}{2}+1.\]
Comme $|z+k| < \frac{1}{2}$ pour au plus un entier naturel $k$, et qu'on a en plus supposé que $|z+k| \geq \frac{1}{8}$ pour tout entier naturel $k$ et $\re(z) < 2-m$, on a
\[\left|\sum_{k=0}^{m-1} \frac{1}{z+k} \right| \leq 10+\sum_{k=0}^{m-3}\frac{1}{-2+(m-k)} < 10+\sum_{j=1}^{m-2} \frac{1}{j} \leq 11+\ln(|z|+1),\]
ce raisonnement valant du moins si $m > 2$ (si $m\leq 2$, l'inégalité vaut même sans le terme logarithmique), pour finalement donner l'inégalité
\[\left|\frac{\Gamma'}{\Gamma}(z) \right|\leq \ln(|z|)+\ln(9)+\ln\left(\sqrt{5}\right)+\frac{\pi}{2}+12,\]
valable pour tout $z$ complexe tel que $|z+k|\geq \frac{1}{8}$ pour tout entier naturel $k$. D'où l'inégalité annoncée.\,$\square$\\

\begin{lem}\label{digamma}Si $\re(s) > -\frac{1}{2}$ et $|s| \geq \frac{1}{8}$, alors \begin{equation}
\left|\frac{\gamma'_\chi}{\gamma_\chi}(s)\right| \leq \frac{n_E}{2}\left(\ln(1+|s|)+{\frac {164}{7}}\right),
\end{equation}\end{lem}

\emph{Démonstration.} Le lemme précédent montre que si $|z|\geq \frac{1}{16}$ et $\re(z) \geq -\frac{1}{4}$, alors:\begin{equation}\left|\frac{\Gamma'}{\Gamma}(z) \right|\leq \ln(|z|)+{\frac {827}{36}}
.\label{inegdigamma}\end{equation}En effet, si on suppose $\re(z)>-\frac{1}{4}$, alors $\re(z+2)>\frac{7}{4}$, donc \[\left|\frac{\Gamma'}{\Gamma}(z+2)\right| \leq \ln(|z|+2)+\frac{\pi}{2}+\frac{4}{7}. \]L'équation fonctionnelle donne alors l'inégalité (\ref{inegdigamma}). En particulier, si \[\gamma_\chi(s) = \left[\pi^{-\frac{s+1}{2}}\Gamma\left(\frac{s+1}{2}\right)\right]^{b(\chi)}\cdot \left[\pi^{-\frac{s}{2}}\Gamma\left(\frac{s}{2}\right)\right]^{a(\chi)},\]
alors 
\[\frac{\gamma'_\chi}{\gamma_\chi}(s) = -\frac{b(\chi)\ln(\pi)}{2} + \frac{b(\chi)}{2}\frac{\Gamma'}{\Gamma}\left(\frac{s+1}{2}\right) - \frac{a(\chi)\ln(\pi)}{2}+\frac{a(\chi)}{2}\frac{\Gamma'}{\Gamma}\left(\frac{s}{2}\right),\]
et le résultat en découle simplement.\,$\square$\\

\emph{Remarque.} Ce même lemme \ref{ineqdigamma} permet de démontrer que si $s=\sigma+i t$ avec $\sigma \geq 2$, alors $\left|\frac{\gamma_\chi'}{\gamma_\chi}(s)\right| \leq \frac{n_E}{2}\left(\ln\left(|t|+\sigma+1\right)+\frac{405}{134}\right)$. Cette inégalité plus précise nous sera utile plus particulièrement dans le prochain lemme (pour $\sigma \geq 1$, on peut remplacer $405$ par $539$, et on l'utilisera dans la démonstration du lemme \ref{sanszero}).

Fort de tous ces résultats, on peut démontrer le premier lemme important de \cite{Lag} sous une forme explicite:

\begin{lem}\label{nchi}Soit $n_\chi(t)$ le nombre de zéros $\rho = \beta+ i \gamma$ de $\mathrm{L}(s,\chi)$ avec $0<\beta<1$ et $|\gamma-t|\leq 1$. Pour tout $t$, on a\[n_\chi(t)+n_\chi(-t) \leq \frac{5}{2}\left[\ln(A(\chi))+n_E\left(\ln(|t|+3)+\frac{1075}{134}\right)\right].\]\end{lem}

\emph{Démonstration.} On part de l'identité (\ref{identiteL'L2}):\[\sum_\rho \left(\frac{1}{s-\rho}+\frac{1}{s-\bar{\rho}}\right) = \frac{\mathrm{L}'}{\mathrm{L}}(s,\chi)+\frac{\mathrm{L}'}{\mathrm{L}}(s,\bar{\chi})+\ln(A(\chi))+2\delta(\chi)\left(\frac{1}{s}+\frac{1}{s-1}\right)+ 2\frac{\gamma'_\chi}{\gamma_\chi}(s).\]
On sait estimer chacune des quantités du membre de droite, d'après ce qui précède. Prenons $s= 2+it$. On a:
\[\left|\sum_\rho \left(\frac{1}{s-\rho}+\frac{1}{s-\bar{\rho}}\right)\right| \leq 3+\ln(A(\chi)) + n_E\left(\ln(|t|+3)+\frac{673}{134}\right),\]puis
\begin{eqnarray}
\sum_\rho \re\left(\frac{1}{s-\rho}+\frac{1}{s-\bar{\rho}}\right) & \geq & \sum_{\substack{\rho \\ |\gamma - t| \leq 1}} \frac{2-\beta}{(2-\beta)^2+(t-\gamma)^2}+\sum_{\substack{\rho \\ |\gamma + t| \leq 1}} \frac{2-\beta}{(2-\beta)^2+(t+\gamma)^2}\nonumber \\
& \geq & \left(\sum_{\substack{\rho \\ |\gamma-t| \leq 1}}+\sum_{\substack{\rho \\ |\gamma+t| \leq 1}}\right) \frac{2}{5} = \frac{2}{5}(n_\chi(t)+n_\chi(-t)).\,\square\nonumber 
\end{eqnarray}

Ce lemme permet d'obtenir d'autres estimations qui informent sur l'importance des zéros des fonctions $\mathrm{L}$ dans la preuve effective du théorème de Chebotarev: des calculs lourds mais sans mystère conduisent à:
\begin{lem}\label{bchi}Pour tout réel $\varepsilon$ tel que $0< \varepsilon \leq 1$, on a\begin{eqnarray}
\left|B(\chi) + \sum_{|\rho| < \varepsilon} \frac{1}{\rho}\right| &\leq & \frac{1}{8}\left(5\pi^2 + 34+\frac{10}{\varepsilon}\right)\ln(A(\chi))\nonumber \\ & \quad +&  \left(\frac {10842}{107}+\frac{1790}{157\varepsilon}\right)n_E \nonumber
\end{eqnarray}\end{lem}
\emph{Démonstration.} On a:
\begin{eqnarray}\left|B(\chi) + \sum_{|\rho| < \varepsilon} \frac{1}{\rho}\right|
   & \leq & \left|B(\chi) + \sum_\rho \left(\frac{1}{2-\rho}+\frac{1}{\rho}\right)\right| \label{sommeB1} \\
   & \quad + & \left|\sum_{|\rho|\geq 1} \left(\frac{1}{2-\rho}+\frac{1}{\rho}\right)\right|  \label{sommeB2}\\
   & \quad + & \left|\sum_{|\rho|<1} \frac{1}{2-\rho}\right| \label{sommeB3}\\
   & \quad + & \left|\sum_{\varepsilon < |\rho| < 1} \frac{1}{\rho}\right|. \label{sommeB4} \end{eqnarray}
   
Les sommes dans (\ref{sommeB2}) et (\ref{sommeB4}) s'estiment de la m\^eme manière, à l'aide de la fonction $n_\chi$ introduite pour le lemme \ref{nchi}. Par exemple, dans le cas de (\ref{sommeB2}), comme $\left|\frac{1}{2-\rho}+\frac{1}{\rho}\right| = \frac{2}{|(2-\rho)\rho|} \leq \frac{2}{|\rho|^2}$, on a
\begin{eqnarray}
\left|\sum_{|\rho|\geq 1} \left(\frac{1}{2-\rho}+\frac{1}{\rho}\right)\right| & \leq & \sum_{\substack{t=-\infty \\ t \textrm{ impair}}}^\infty \sum_{\substack{|\rho|\geq 1 \\ t \leq \gamma \leq t+2}} \frac{2}{|\rho|^2} \nonumber \\
 & \leq & 2\sum_{j=0}^\infty \frac{n_\chi(2j+2)+n_\chi(-(2j+2))}{(2j+1)^2}+2n_\chi(0) \nonumber\\
 & \leq & 5\left(\sum_{j=0}^\infty \frac{1}{(2j+1)^2}+\frac{1}{2}\right)\left(\ln(A(\chi))+\frac{1075}{134}n_E\right)\nonumber\\
 & \quad + &5n_E\left(\ln(5)+\int_0^\infty \frac{\ln(2t+5)}{(2t+1)^2}\ud t+\frac{1}{2}\ln(3)\right)\nonumber\\
 & = & \frac{5}{8}(\pi^2+4)\ln(A(\chi)) \nonumber\\
 & \quad +& 5\left[\frac{1075}{134\cdot 8}\left(\pi^2+4\right)+\left(\frac{13}{8}\ln(5)+\frac{\ln(3)}{2}\right)\right]n_E. \nonumber
\end{eqnarray}
Semblablement, on trouve \[\left|\sum_{\varepsilon < |\rho| < 1} \frac{1}{\rho}\right|\leq \frac{n_\chi(0)}{\varepsilon}\leq \frac{5}{4\varepsilon}\left[\ln(A(\chi))+n_E\left(\ln(3)+\frac{1075}{134}\right)\right].\]
Si $|\rho|<1$, alors $|2-\rho|>1$, donc \[\left|\sum_{|\rho|<1} \frac{1}{2-\rho}\right| \leq \frac{5}{4}\left[\ln(A(\chi))+n_E\left(\ln(3)+\frac{1075}{134}\right)\right],\]et il ne reste plus que la somme (\ref{sommeB1}) à évaluer. Pour cela, on pose $s=2$ dans (\ref{identiteL'L}), et on estime tous les termes qui nous intéressent grâce aux lemmes \ref{L'surL} et \ref{digamma} (plus précisément, la remarque qui suit le lemme \ref{digamma}):
\[ \left|B(\chi) + \sum_\rho \left(\frac{1}{2-\rho}+\frac{1}{\rho}\right)\right| \leq \frac{3}{2}+\frac{n_E}{2}\left(\ln(3)+\frac{673}{134}\right) + \frac{1}{2}\ln(A(\chi)).\]
En regroupant toutes ces estimations, on obtient le lemme annoncé
.\,$\square$\\

Dans le même esprit de preuve:
\begin{lem}\label{L'/L}Si $s = \sigma+it$ avec $-\frac{1}{2}\leq \sigma \leq 3$ et $|s|\geq \frac{1}{8}$, alors:
\begin{eqnarray}
\left|\frac{\mathrm{L}'}{\mathrm{L}}(s,\chi)+\frac{\delta(\chi)}{s-1}-\sum_{\substack{\rho\\ |\gamma - t| \leq 1}} \frac{1}{s-\rho}\right| & \leq & \frac{5}{4}\left(1+\frac{7}{4}\pi^2\right)\ln(A(\chi)) \nonumber \\
 &\quad + & \frac{n_E}{2}\ln(|t|+5)\left(57+\frac{35}{|t|+4}\right)\nonumber \\
 &\quad + & \frac{50096}{255}n_E + \frac{53}{6}. \nonumber
\end{eqnarray}
\end{lem}

\emph{Démonstration.} On évalue l'identité (\ref{identiteL'L}) en $\sigma + it$ puis en $3+it$, et on soustrait les deux égalités résultantes. Ainsi, $B(\chi)$ est éliminé, et on a:
\begin{eqnarray}
\frac{\mathrm{L}'}{\mathrm{L}}(s,\chi)-\frac{\mathrm{L}'}{\mathrm{L}}(3+it,\chi) & = & \sum_\rho \left(\frac{1}{s-\rho}-\frac{1}{3+it-\rho}\right)-\frac{\gamma_\chi'}{\gamma_\chi}(s)+\frac{\gamma_\chi'}{\gamma_\chi}(3+it)\nonumber \\
& \quad -& \delta(\chi)\left(\frac{1}{s}+\frac{1}{s-1}-\frac{1}{2+it}-\frac{1}{3+it}\right).\nonumber
\end{eqnarray}
La seule quantité à réellement poser problème, dans la majoration qui nous intéresse, est la somme sur les zéros $\rho$. L'application des lemmes \ref{L'surL} et \ref{digamma} permet de majorer de la façon suivante:
\begin{eqnarray}
\left|\frac{\mathrm{L}'}{\mathrm{L}}(s,\chi)+\frac{\delta(\chi)}{s-1}-\sum_{\substack{\rho\\ |\gamma - t| \leq 1}} \frac{1}{s-\rho}\right| & \leq & \sum_{\substack{\rho \\ |\gamma - t| > 1}} \left|\frac{1}{s-\rho}-\frac{1}{3+it-\rho}\right| \label{sommelemme561}
\\ & \quad + &  \sum_{\substack{\rho \\ |\gamma - t| \leq 1}} \left|\frac{1}{3+it-\rho}\right| \label{sommelemme562}
\\ & \quad + & \frac{53}{6} + \left(\ln(|t|+4)+\frac{24811}{1876}\right)n_E.\nonumber
\end{eqnarray}
On a $|3+it-\rho| > 2$ pour tout zéro $\rho$ non trivial (car $0<\re(\rho)<1$), donc \[\sum_{\substack{\rho \\ |\gamma-t|\leq 1}}\left|\frac{1}{3+it-\rho}\right| \leq \frac{1}{2}n_\chi(t), \]et on utilise alors le lemme \ref{nchi} pour majorer (\ref{sommelemme562}). \`A présent, (\ref{sommelemme561}) se calcule progressivement (on note $s = \sigma + it$ et $\rho = \beta+i\gamma$):
\begin{eqnarray}
\sum_{\substack{\rho \\ |\gamma - t| > 1}} \left|\frac{1}{s-\rho}-\frac{1}{3+it-\rho}\right|  & = & \sum_{\substack{\rho \\ |\gamma - t| > 1}} \frac{3-\sigma}{|s-\rho||3+it-\rho|}  \nonumber\\
                    & \leq & \frac{7}{2}\sum_{k=1}^\infty \sum_{\substack{\rho \\ |\gamma-(t+2k)| \leq 1}}\frac{1}{|t-\gamma|^2} \nonumber \\
                    & \quad  +& \frac{7}{2}\sum_{k=1}^{\infty} \sum_{\substack{\rho \\ |\gamma-(t-2k)| \leq 1}}\frac{1}{|t-\gamma|^2} \nonumber \\
                    & \leq & \frac{7}{2}\sum_{k=1}^\infty \frac{n_\chi(t+2k)+n_\chi(t-2k)}{(2k-1)^2} \nonumber \\
                    & \leq & \frac{35\pi^2}{16}\left(\ln(A(\chi)) + \frac{1075}{134}n_E\right) \nonumber \\
                    & \quad + & \frac{35}{2}\sum_{k=1}^\infty \frac{\ln(|t|+2k+3)}{(2k-1)^2}n_E. \nonumber
\end{eqnarray}
Cette dernière somme se majore, par exemple, par une comparaison série-intégrale; on obtient ainsi $\frac{\ln(|t|+5)}{2}\left(3+\frac{1}{|t|+4}\right)$ comme majorant. On en déduit le lemme, en regroupant les différentes estimations faites ici. Par commodité, on utilise la majoration $\frac{24811}{1876}+\frac{1075}{134}\left(\frac{5}{4}+\frac{35\pi^2}{16}\right)\leq \frac {50096}{255}$.\,$\square$

\section{L'intégrale sur un contour}

Maintenant, on se charge d'évaluer $I_C(x,T)$ en passant par l'évaluation de \begin{equation}\label{defIchi}
I_\chi(x,T) = \frac{1}{2\pi i}\int_{\sigma_0-iT}^{\sigma_0+iT} \frac{x^s}{s}\frac{\mathrm{L}'}{\mathrm{L}}(s,\chi)\ud s
\end{equation}pour chaque caractère $\chi$ de $H$. À présent, on a besoin que $T$ soit différent des zéros de toutes les fonctions $\mathrm{L}(s,\chi)$ (cette condition s'éclairera plus tard). On introduit un nouveau paramètre $U$, vérifiant $U = j+\frac{1}{2}$ pour un certain entier naturel $j$ (plus tard, on fera tendre $U$ vers $+\infty$), et on définit
\begin{equation}I_\chi(x,T,U) = \frac{1}{2\pi i}\int_{B_{T,U}} \frac{x^s}{s}\frac{\mathrm{L}'}{\mathrm{L}}(s,\chi)\ud s,\label{defichixtu}\end{equation}
où $B_{T,U}$ est le rectangle (orienté dans le sens trigonométrique) dont les sommets sont $\sigma_0-iT$, $\sigma_0+iT$, $-U+iT$ et $-U-iT$. Cette intégrale s'exprime simplement à l'aide des singularités de l'intégrande, mais ceci attendra encore un peu. Dans cette section, on montrera que $R_\chi(x,T,U) = I_\chi(x,T,U) - I_\chi(x,T)$ est petit. On divise $R_\chi(x,T,U)$ en trois intégrales, celle verticale
\begin{equation}\label{defV}
V_\chi(x,T,U) = \frac{1}{2\pi}\int_T^{-T} \frac{x^{-U+iT}}{-U+iT} \frac{\mathrm{L}'}{\mathrm{L}}(-U+it,\chi)\ud t,
\end{equation}
et les deux intégrales horizontales
\begin{equation}\label{defH}
H_\chi(x,T,U) = \frac{1}{2\pi i}\int_{-U}^{-1/4} \left(\frac{x^{\sigma-iT}}{\sigma-iT}\frac{\mathrm{L}'}{\mathrm{L}}(\sigma-iT,\chi) - \frac{x^{\sigma+iT}}{\sigma+iT} \frac{\mathrm{L}'}{\mathrm{L}}(\sigma+iT,\chi)\right)\ud \sigma
\end{equation}
et
\begin{equation}\label{defHstar}
H_\chi^*(x,T) = \frac{1}{2\pi i}\int_{-1/4}^{\sigma_0} \left(\frac{x^{\sigma-iT}}{\sigma-iT}\frac{\mathrm{L}'}{\mathrm{L}}(\sigma-iT,\chi) - \frac{x^{\sigma+iT}}{\sigma+iT} \frac{\mathrm{L}'}{\mathrm{L}}(\sigma+iT,\chi)\right)\ud \sigma
\end{equation}
Pour borner $\frac{\mathrm{L}'}{\mathrm{L}}$, et donc estimer ces intégrales, on a besoin du lemme \ref{ineqdigamma}, qui renseigne sur la taille de $\frac{\Gamma'}{\Gamma}$.


On en déduit:
\begin{lem}Si $s = \sigma + it$ avec $\sigma \leq -\frac{1}{4}$ et $|s+m| \geq \frac{1}{4}$ pour tout entier naturel $m$, alors
\[\left|\frac{\mathrm{L}'}{\mathrm{L}}(s,\chi)\right| \leq \ln(A(\chi)) +n_E\left(\ln(|s|+2)+\frac{19683}{812}\right).\]\end{lem}

\emph{Démonstration.} Partant de l'équation fonctionnelle vérifiée par $\mathrm{L}(\cdot,\chi)$, on a\[\frac{\mathrm{L}'}{\mathrm{L}}(s,\chi) = -\frac{\mathrm{L}'}{\mathrm{L}}(1-s,\bar{\chi}) -\ln(A(\chi)) - \frac{\gamma_\chi'}{\gamma_\chi}(1-s) - \frac{\gamma_\chi'}{\gamma_\chi}(s),\]
et $\re(1-s) \geq \frac{5}{4}$, ce qui permet d'utiliser le lemme \ref{L'surL} pour borner $\frac{\mathrm{L}'}{\mathrm{L}}(1-s,\bar{\chi})$. De plus, le lemme \ref{ineqdigamma} donne une autre version du lemme \ref{digamma}, pour $z$ tel que $|z+k| \geq \frac{1}{8}$ pour tout entier naturel $k$: on obtient l'inégalité similaire $\left|\frac{\gamma_\chi'}{\gamma_\chi}(z)\right| \leq \frac{n_E}{2}\left(\ln(1+|z|)+\frac{989}{58}\right)$. Alors,
\[\left|\frac{\mathrm{L}'}{\mathrm{L}}(s,\chi)\right| \leq \ln(A(\chi)) - \frac{n_E}{\sigma} + n_E\left(\ln(|s|+2)+\frac{16435}{812}\right),\]
d'où le résultat.\,$\square$\\

Reprenons: les intégrales verticale (\ref{defV}) et horizontale (\ref{defH}) sont estimées à l'aide de l'inégalité ci-dessus. On a, pour $x\geq 2$, $T \geq 1$ et $U \geq \frac{1}{2}$,

\begin{eqnarray}|V_\chi(x,T,U)|& \leq & \frac{x^{-U}}{2\pi U}\int_{-T}^T \left|\frac{\mathrm{L}'}{\mathrm{L}}(-U+it,\chi)\right|\ud t \nonumber \\
& \leq & \frac{x^{-U} T}{\pi U}\left(\ln(A(\chi)) +n_E\left(\ln(U+T+2)+\frac{19683}{812}\right)\right),\label{vchi}\end{eqnarray}
et, sous ces mêmes conditions,
\begin{eqnarray}|H_\chi(x,T,U)| &\leq & \frac{1}{\pi T}\int_{-\infty}^{-\frac{1}{4}} x^\sigma\left(\ln(A(\chi)) +n_E\left(\ln(|\sigma|+T+2)+\frac{19683}{812}\right)\right)\ud \sigma \nonumber \\
& \leq & \frac{x^{-1/4}}{\pi T \ln(x)} \left[\ln(A(\chi))+\left(\ln\left(T+\frac{9}{4}\right)+\frac{19683}{812}\right)n_E\right] \nonumber \\
& \quad + & \frac{n_E x^{-1/4}}{\pi T (\ln(x))^2}\frac{1}{\left(T+\frac{9}{4}\right)}, \label{hchi}\end{eqnarray}
Estimer l'intégrale $H_\chi^*(x,T)$ est un peu plus ardu, et nécessite le lemme \ref{L'/L}. Il nous indique que
\begin{eqnarray}
\left|\frac{\mathrm{L}'}{\mathrm{L}}(\sigma+iT,\chi)-\sum_{\substack{\rho \\ |\gamma-T| \leq 1}} \frac{1}{\sigma+iT-\rho} \right| &\leq & \frac {571}{25}\ln(A(\chi)) \nonumber \\
 &\quad + & n_E\left(\frac{57}{2}\ln(T+5)+\frac{5921}{28}\right) \nonumber \end{eqnarray} pour $\sigma$ dans $\left[-\frac{1}{4},\sigma_0\right]$, $x \geq 2$ et $T\geq 2$. On a la même estimation en $\sigma-iT$. Alors,
\begin{eqnarray}
\Big|H_\chi^*(x,T) & - & \frac{1}{2\pi i }\int_{-1/4}^{\sigma_0} \Big(\frac{x^{\sigma-iT}}{\sigma-iT}\sum_{\substack{\rho \\ |\gamma + T| \leq 1}} \frac{1}{\sigma - iT - \rho} \nonumber \\
                   & \qquad - & \frac{x^{\sigma+iT}}{\sigma+iT}\sum_{\substack{\rho \\ |\gamma-T| \leq 1}} \frac{1}{\sigma + iT - \rho}\Big)\ud \sigma\Big| \nonumber \\
& \leq & \frac{ex-x^{-1/4}}{\pi T\ln(x)}\left[\frac{571}{25}\ln(A(\chi))+n_E\left(\frac{57}{2}\ln(T+5)+\frac{5921}{28}\right)\right].\nonumber
\end{eqnarray}
Pour en finir avec l'estimation de $H_\chi^*(x,T)$, on doit encore vérifier que l'intégrale ci-dessus reste assez «~petite~».

\begin{lem}Soit $\rho = \beta + i\gamma$, avec $0<\beta<1$ et $\gamma \neq t$. Si $|t|\geq 2$, $x \geq 2$ et $1 < \sigma_1 \leq 3$, alors: \[\left|\int_{-1/4}^{\sigma_1} \frac{x^{\sigma+it}}{(\sigma+it)(\sigma+it-\rho)}\ud \sigma\right| \leq \left(\sigma_1+\frac{9}{4}\right)\frac{x^{\sigma_1}}{(|t|-1)(\sigma_1-\beta)}.\]\end{lem}

\emph{Démonstration.} Supposons d'abord que $\gamma > t$. Soit $B$ le rectangle (orienté dans le sens trigonométrique) dont les sommets ont pour affixes $\sigma_1+i(t-1)$, $\sigma_1+it$, $-\frac{1}{4}+it$ et $-\frac{1}{4}+i(t-1)$. Le théorème de Cauchy assure que\[\int_B \frac{x^s}{s(s-\rho)}\ud s = 0,\]puisque l'intégrande n'a pas de singularité à l'intérieur du rectangle. En outre, sur les arêtes du rectangle, à l'exception de celle joignant $-\frac{1}{4}+it$ et $\sigma_1+it$, on peut le majorer par $\frac{x^{\sigma_1}}{(|t|-1)(\sigma_1-\beta)}$. D'où le résultat pour $\gamma > t$. On procède de même si $\gamma < t$, en changeant $i(t-1)$ en $i(t+1)$ dans les affixes des sommets du rectangle.\,$\square$\\

Ceci prouve que \begin{eqnarray} \Big| \frac{1}{2\pi i} \int_{-1/4}^{\sigma_0}\frac{x^{\sigma-iT}}{\sigma-iT}&\Big( &\sum_{\substack{\rho \\ |\gamma+T| \leq 1}} \frac{1}{\sigma-iT-\rho}\Big)\ud  \sigma \Big| \nonumber \\
 & \leq & \frac{\sigma_0+\frac{9}{4}}{2\pi}n_\chi(-T)\frac{x^{\sigma_0}}{(T-1)(\sigma_0-1)} \nonumber \\
& \leq & \frac{5\left(13\ln(x)+4\right)ex}{16\pi (T-1)}\left[\ln(A(\chi))+n_E\left(\ln(T+3)+\frac{1075}{134}\right)\right]\nonumber\end{eqnarray}pour $x \geq 2$ et $T\geq 2$. On a le même résultat pour les zéros indicés grâce à la condition $|\gamma - T|\leq 1$. Si on suppose l'hypothèse de Riemann généralisée, le terme en $\ln(x)$ disparaît, d'après \cite{Lag}, page 445. On ne s'en chargera pas ici.\\
En fin de compte, on obtient
\begin{eqnarray}|H_\chi^* (x,T)| & \leq &  \left(\frac{5\left(13\ln(x)+4\right)ex}{8 (T-1)}+\frac{571}{25}\frac{ex-x^{-1/4}}{T\ln(x)}\right)\frac{\ln(A(\chi))}{\pi}\nonumber\\
& \quad +& \frac{n_E}{\pi}\Big[\left(\frac{57}{2}\frac{ex-x^{-1/4}}{T\ln(x)}+\frac{5\left(13\ln(x)+4\right)ex}{8(T-1)}\right)\ln(T+5)\nonumber\\
& \quad +& \frac{5375\left(13\ln(x)+4\right)ex}{2144 (T-1)}+\frac{5921}{28}\frac{ex-x^{-1/4}}{ T\ln(x)}\Big]\label{hchistar}\end{eqnarray}

En combinant (\ref{vchi}), (\ref{hchi}) et (\ref{hchistar}), on obtient le résultat principal de la section, à savoir,

\begin{eqnarray}
|I_\chi(x,T) - I_\chi(x,T,U)| &\leq & \frac{65e}{8\pi}\frac{x\ln(x)}{T-1}\left[\ln(A(\chi))+n_E\left(\ln(T+5)+\frac{1075}{268}\right)\right]
\nonumber \\
& \quad + & \frac{5e}{2\pi}\frac{x}{T-1}\left[\ln(A(\chi))+n_E\left(\ln(T+5)+\frac{1075}{268}\right)\right]\nonumber\\
& \quad + & \frac{e}{\pi}\frac{x}{T\ln(x)}\left[\frac{571}{25}\ln(A(\chi))+n_E \left(\frac{57}{2}\ln(T+5)+\frac{5921}{28}\right)\right]\nonumber\\
&\quad  + &\frac{x^{-U} T}{\pi U}\left[\ln(A(\chi)) +n_E\left(\ln(U+T+2)+\frac{19683}{812}\right)\right]\nonumber\\
& \quad + & \frac{4n_E x^{-1/4}}{17\pi T(\ln(x))^2}.\label{troncatureichi}
\end{eqnarray}
\section{La formule explicite}

On en arrive enfin à une formule explicite pour $\psi_C$ en fonction des zéros $\rho$. On revient à la définition de $I_\chi(x,T,U)$ donnée en (\ref{defichixtu}), où on rappelle que $x \geq 2$ et $U = j+\frac{1}{2}$ pour un certain entier naturel $j$. Soit $T \geq 2$ différent de l'ordonnée de tout zéro de tout $\mathrm{L}(s,\chi)$. Par le théorème de Cauchy, $I_\chi(x,T,U)$ égale la somme des résidus de l'intégrande aux pôles à l'intérieur de $B_{T,U}$. Pour un décompte détaillé des résidus, on renvoie à \cite{Lag}, pages 446--448: on a en tout cas
\begin{eqnarray}
I_\chi(x,T,U) & = & -\delta(\chi)x + \sum_{\substack{\rho \\ |\gamma|<T}} \frac{x^\rho}{\rho} - b(\chi)\sum_{1\leq m \leq \frac{U+1}{2}} \frac{x^{1-2m}}{2m-1}\nonumber \\
& \quad - & a(\chi)\sum_{1\leq m \leq \frac{U}{2}} \frac{x^{-2m}}{2m} + r(\chi) + (a(\chi)-\delta(\chi))\ln(x),\label{passagelimite}
\end{eqnarray}
où\[r(\chi) = B(\chi) - \frac{1}{2}\ln(A(\chi))+\frac{n_E}{2}\ln(\pi) + \delta(\chi) - \frac{b(\chi)}{2}\frac{\Gamma'}{\Gamma}\left(\frac{1}{2}\right) - \frac{a(\chi)}{2}\frac{\Gamma'}{\Gamma}(1).\]
Quand $U \to +\infty$, (\ref{troncatureichi}) et (\ref{passagelimite}) donnent la formule explicite suivante, valide pour tout $x \geq 2$ et tout $T \geq 2$ différent de l'ordonnée d'un zéro.
\begin{eqnarray}
\Big|I_\chi(x,T) + \delta(\chi)x & - & \sum_{\substack{\rho \\ |\gamma|<T}} \frac{x^\rho}{\rho}-r(\chi)-(a(\chi)-\delta(\chi))\ln(x) \nonumber\\
& \quad - & \frac{n_E}{2}\ln\left(1-\frac{1}{x}\right)+\frac{1}{2}(b(\chi)-a(\chi))\ln\left(1+\frac{1}{x}\right)\Big|\nonumber\\
& \leq &\frac{65e}{8\pi}\frac{x\ln(x)}{T-1}\left[\ln(A(\chi))+n_E\left(\ln(T+5)+\frac{1075}{268}\right)\right]
\nonumber \\
& \quad + & \frac{5e}{2\pi}\frac{x}{T-1}\left[\ln(A(\chi))+n_E\left(\ln(T+5)+\frac{1075}{268}\right)\right]\nonumber\\
& \quad + & \frac{e}{\pi}\frac{x}{T\ln(x)}\left[\frac{571}{25}\ln(A(\chi))+n_E \left(\frac{57}{2}\ln(T+5)+\frac{5921}{28}\right)\right]\nonumber\\
& \quad + & \frac{4n_E x^{-1/4}}{17\pi T(\ln(x))^2}.
\end{eqnarray}
On a alors presque ce qu'on veut:
\begin{theo}\label{formulexplicite}Si $x \geq 2$ et $T\geq 2$, alors \begin{eqnarray}
\left|\psi_C(x) - \frac{|C|}{|G|}x +S(x,T)\right| & \leq & \frac{|C|}{|G|}\Big(\frac{145}{6}n_L \frac{x(\ln(x))^2}{T}\nonumber\\
& \quad + & \frac{65e}{8\pi}\frac{x\ln(x)}{T-1}\left[\ln(d_L)+n_L\left(\ln(T+5)+\frac{1075}{268}\right)\right]
\nonumber \\
& \quad + & \frac{5}{2}\left(\frac{e}{\pi}+1\right)\frac{x}{T-1}\left[\ln(d_L)+n_L\left(\ln(T+5)+\frac{1075}{268}\right)\right]\nonumber\\
& \quad + & \frac{e}{\pi}\frac{x}{T\ln(x)}\left[\frac{571}{25}\ln(d_L)+n_L \left(\frac{57}{2}\ln(T+5)+\frac{5921}{28}\right)\right]\nonumber\\
& \quad + & \ln(x)\left[\frac{2}{\ln(2)}\ln(d_L)+\frac{77}{4}n_L\right] + \frac{94}{7}\ln(d_L) + \frac{3817}{30}n_L\nonumber\\
& \quad + & \frac{4n_L x^{-1/4}}{17\pi T(\ln(x))^2}\Big),\nonumber
\end{eqnarray}où \begin{equation}
S(x,T) = \frac{|C|}{|G|}\sum_\chi \bar{\chi}(g) \left(\sum_{\substack{\rho\\ |\gamma| < T}} \frac{x^\rho}{\rho} - \sum_{\substack{\rho\\ |\rho| < \frac{1}{2}}} \frac{1}{\rho}\right),
\end{equation}les deux dernières sommes étant sur les zéros non triviaux de $\mathrm{L}(\cdot,\chi)$.
\end{theo}

\emph{Démonstration.} D'après le lemme \ref{bchi}, on a \[\left|r(\chi)+\sum_{|\rho| < \frac{1}{2}}\frac{1}{\rho}\right| \leq \left(\frac{5\pi^2}{8}+\frac{29}{4}\right)\ln(A(\chi)) + \left(\frac {2102053}{16799}+\frac{1}{2}\ln(\pi)+\frac{1}{2}\gamma_0+\ln(2)\right)n_E,\]
donc \begin{eqnarray}\Big|I_\chi(x,T) + \delta(\chi)x & - & \sum_{\substack{\rho \\|\gamma|<T}} \frac{x^\rho}{\rho}  - \sum_{|\rho|<\frac{1}{2}} \frac{1}{\rho}\Big| \nonumber \\ &\leq &  \frac{65e}{8\pi}\frac{x\ln(x)}{T-1}\left[\ln(A(\chi))+n_E\left(\ln(T+5)+\frac{1075}{268}\right)\right]
\nonumber \\
& \quad + & \frac{5e}{2\pi}\frac{x}{T-1}\left[\ln(A(\chi))+n_E\left(\ln(T+5)+\frac{1075}{268}\right)\right]\nonumber\\
& \quad + & \frac{e}{\pi}\frac{x}{T\ln(x)}\left[\frac{571}{25}\ln(A(\chi))+n_E \left(\frac{57}{2}\ln(T+5)+\frac{5921}{28}\right)\right]\nonumber\\
& \quad + & n_E \ln(x) + \frac{94}{7}\ln(A(\chi)) + \frac{3817}{30}n_E\nonumber\\
& \quad + & \frac{4n_E x^{-1/4}}{17\pi T(\ln(x))^2}. \nonumber
\end{eqnarray}
Alors, (\ref{rappelIC}) et (\ref{defIchi}) donnent, pour $T \neq \mathrm{Im}(\rho)$ pour tout zéro $\rho$ de tout $\mathrm{L}(s,\chi)$, la même majoration (à une multiplication par $\frac{|C|}{|G|}$ près) pour \[\left|I_C(x,T)- \frac{|C|}{|G|}\sum_\chi \bar{\chi}(g)\left(\delta(\chi)x - \sum_{\substack{\rho \\ |\gamma|<T}} \frac{x^\rho}{\rho} - \sum_{|\rho|<\frac{1}{2}} \frac{1}{\rho}\right)\right|\]en remplaçant $\ln(A(\chi))$ par $\ln(d_L)$ et $n_E$ par $n_L$, en vertu des formules \[\sum_\chi \ln(A(\chi)) = \ln(d_L)\textrm{ (formule du discriminant) et }n_E \cdot [L:E] = n_L. \]Comme $\psi_C(x,T) = I_C(x,T)+R_1(x,T)$,

\begin{eqnarray}
\left|\psi_C(x) - \frac{|C|}{|G|}x +S(x,T)\right|
& \leq & \frac{|C|}{|G|}\Big(\frac{145}{6}n_L \frac{x(\ln(x))^2}{T}\nonumber\\
& \quad + & \frac{65e}{8\pi}\frac{x\ln(x)}{T-1}\left[\ln(d_L)+n_L\left(\ln(T+5)+\frac{1075}{268}\right)\right]
\nonumber \\
& \quad + & \frac{5e}{2\pi}\frac{x}{T-1}\left[\ln(d_L)+n_L\left(\ln(T+5)+\frac{1075}{268}\right)\right]\nonumber\\
& \quad + & \frac{e}{\pi}\frac{x}{T\ln(x)}\left[\frac{571}{25}\ln(d_L)+n_L \left(\frac{57}{2}\ln(T+5)+\frac{5921}{28}\right)\right]\nonumber\\
& \quad + & \ln(x)\left[\frac{2}{\ln(2)}\ln(d_L)+\frac{77}{4}n_L\right] + \frac{94}{7}\ln(d_L) + \frac{3817}{30}n_L\nonumber\\
& \quad + & \frac{4n_L x^{-1/4}}{17\pi T(\ln(x))^2}\Big),\nonumber
\end{eqnarray}
ce qui donne le théorème si $T \neq \mathrm{Im}(\rho)$ pour tout zéro $\rho$ de fonction $\mathrm{L}$. Maintenant, si $T = \mathrm{Im}(\rho_T)$ pour un certain zéro $\rho_T$ (ou plusieurs; il y en a au plus $\sum\limits_\chi n_\chi(T)$), on a, pour $\varepsilon$ assez proche de 0, 
\[S(x,T+\varepsilon) = S(x,T)+\frac{|C|}{|G|}\sum_\chi \bar{\chi}(g)\sum_{\substack{\rho_T \\ T = \mathrm{Im}(\rho_T)}}\frac{x^{\rho_T}}{\rho_T}\]
et $T+\varepsilon$ différent de tout zéro de fonction $\mathrm{L}$. Alors,

\begin{eqnarray}
\left|\psi_C(x) - \frac{|C|}{|G|}x +S(x,T)\right| & \leq & \left|\psi_C(x) - \frac{|C|}{|G|}x +S(x,T+\varepsilon)\right|\nonumber \\
& \quad + & \frac{|C|}{|G|}\sum_\chi n_\chi(T) \frac{x}{T} \nonumber
\end{eqnarray}
On sait estimer la dernière somme grâce au lemme \ref{nchi} et la formule du discriminant, tandis que la première quantité du membre de droite s'évalue à l'aide des calculs qui précèdent. On obtient la majoration du théorème en faisant tendre $\varepsilon$ vers 0.\,$\square$\\

C'est le résultat principal du papier: on a exprimé $\psi_C(x)$ en fonction d'un terme principal $\frac{|C|}{|G|}x$, de $S(x,T)$ et un terme d'erreur relativement petit. Il reste à estimer $S(x,T)$. Si on suppose la justesse de l'hypothèse de Riemann généralisée, on peut avoir une bonne borne à partir des résultats déjà établis; si on veut un résultat inconditionnel, on doit montrer que les zéros $\rho$ ne s'approchent pas trop de la droite $\re(s)=1$.
\section{Les régions sans zéros}

Cette section concerne les zéros exceptionnels que les fonctions $\zeta_L$ peuvent éventuellement avoir. Pour rappel, comme $\zeta_L = \prod_\chi \mathrm{L}(\cdot,\chi)$, un zéro pour $\zeta_L$ en entraîne un pour au moins une des fonctions $\mathrm{L}(\cdot,\chi)$; dans le cas du lemme \ref{presquesanszero} à venir, il n'y a en vérité qu'une seule fonction $\mathrm{L}(\cdot,\chi)$ éventuellement concernée. Seul un résultat nécessite d'être explicité, c'est le lemme 8.1 de \cite{Lag}, que voici:

\begin{lem}\label{sanszero}La fonction $\zeta_L$ n'a pas de zéros $\rho = \beta + i\gamma$ dans la région délimitée par les conditions \[|\gamma| \geq \frac{1}{1+4\ln(d_L)}\]\begin{center}et\end{center}\[\beta \geq 1 - (7-4\sqrt{3})\left(22\ln(d_L)+n_L\left(\frac{5}{2}\ln(|\gamma|+3)+\frac{1078}{67}+2\ln(3)\right)+\frac{15}{2}\right)^{-1}.\]\end{lem}

\emph{Démonstration.} Partant de la formule du produit eulérien pour $\zeta_L$, on a pour $\sigma = \re(s)>1$,
\[-\frac{\zeta_L'}{\zeta_L}(s) = \sum_{m=1}^\infty \Lambda_L(m)m^{-s},\]
où
\[\Lambda_L(m) = \left\{\begin{array}{l}\left(\sum\limits_{\substack{\Pre | p \\ \mathfrak{f}(\Pre/p) | r}} \mathfrak{f}(\Pre/p)\right) \ln(p)\textrm{ si }n = p^r, p\textrm{ premier, }r \geq 1,\\0 \textrm{ sinon.}\end{array}\right. \]
En particulier, $\Lambda_L(m)\geq 0$ pour tout entier naturel non nul $m$. Par conséquent, en vertu de l'identité traditionnelle 
\[\forall \theta \in \R, \quad 3+4\cos(\theta)+\cos(2\theta) = 2(1+\cos(\theta))^2,\]
on sait que
\begin{eqnarray}\re\left(-3\frac{\zeta_L'}{\zeta_L}(\sigma) - 4\frac{\zeta_L'}{\zeta_L}(\sigma+it)-\frac{\zeta_L'}{\zeta_L}(\sigma+2it)\right) &=& 2\sum_{m=1}^\infty \Lambda_L(m)m^{-\sigma}(1+\cos(t\ln(m)))^2 \nonumber \\
& \geq & 0. \nonumber\end{eqnarray}
Or, si on prend pour $\chi$ le caractère trivial, l'égalité (\ref{identiteL'L}) montre que 
\[2\frac{\zeta_L'}{\zeta_L}(s) = \sum_\rho \left(\frac{1}{s-\rho}+\frac{1}{s-\bar{\rho}}\right)-\ln(d_L)-2\left(\frac{1}{s}+\frac{1}{s-1}\right)-2\frac{\gamma_L'}{\gamma_L}(s),\]
où la somme est indicée par les zéros $\rho$ de $\zeta_L$. Si $\re(s)>1$, alors $\re\left(\frac{1}{s-\rho}\right)>0$ pour tout $\rho$. \`A présent, considérons $\beta+i\gamma$ un zéro tel que $|\gamma|\geq \frac{1}{1+4\ln(d_L)}$. Alors, pour $2 \geq \sigma>1$, d'après l'inégalité citée dans la remarque qui précède le lemme \ref{nchi},
\begin{eqnarray}-\frac{\zeta_L'}{\zeta_L}(\sigma) &\leq & \frac{1}{\sigma-1}+\frac{1}{\sigma} + \frac{1}{2}\ln(d_L)+\frac{\gamma_L'}{\gamma_L}(\sigma), \nonumber \\
& \leq & \frac{1}{\sigma-1}+ 1+\frac{1}{2}\ln(d_L)+ \frac{n_L}{2}\left(\ln(3)+\frac{539}{134}\right). \nonumber \end{eqnarray}
De même,
\begin{eqnarray}
-\re\left(\frac{\zeta'_L}{\zeta_L}(\sigma+2i\gamma)\right) & \leq & \frac{1}{2}\ln(d_L) + \re\left(\frac{1}{\sigma+2i\gamma-1}+\frac{1}{\sigma+2i\gamma}\right)+\re\left(\frac{\gamma_L'}{\gamma_L}(\sigma+2i\gamma)\right) \nonumber \\
& \leq & \frac{5}{2}\ln(d_L)+\frac{n_L}{2}\left(\ln(2|\gamma|+3)+\frac{539}{134}\right)+\frac{1}{2}, \nonumber
\end{eqnarray}
et, en conservant la contribution du zéro $\rho = \beta + i\gamma$,
\[-\re\left(\frac{\zeta'_L}{\zeta_L}(\sigma+i\gamma)\right) \leq \frac{9}{2}\ln(d_L)+\frac{n_L}{2}\left(\ln(|\gamma|+3)+\frac{539}{134}\right)+1-\frac{1}{\sigma-\beta}.\]
Ces inégalités, et le fait souligné ci-dessus que tous les coefficients de Dirichlet de $-\frac{\zeta_L'}{\zeta_L}$ soient positifs, impliquent que pour tout $\sigma>1$,
\[\frac{4}{\sigma-\beta} < \frac{3}{\sigma-1}+22\ln(d_L)+n_L\left(\frac{5}{2}\ln(|\gamma|+3)+4\cdot\frac{539}{134}+2\ln(3)\right)+\frac{15}{2},\]
et il suffit alors de prendre \[\sigma = 1+(2\sqrt{3}-3)\left(22\ln(d_L)+n_L\left(\frac{5}{2}\ln(|\gamma|+3)+4\cdot\frac{539}{134}+2\ln(3)\right)+\frac{15}{2}\right)^{-1} \]pour obtenir le résultat du lemme.\,$\square$

\begin{lem}\label{presquesanszero}Si $n_L > 1$, alors $\zeta_L$ a au plus un zéro $\rho = \beta + i\gamma$ dans la région délimitée par les conditions \[|\gamma| \leq \frac{1}{4\ln(d_L)}\textrm{ et }\beta \geq 1 - \frac{1}{4\ln(d_L)}.\]Ce zéro, s'il existe, est réel et simple, et correspond à un unique caractère (réel) $\chi_0$ qui annule $\mathrm{L}(\cdot,\chi_0)$. Par ailleurs, si $n_L=1$, alors $\zeta_L$ n'admet pas de zéro tel que $|\gamma|<14$.\end{lem}

\emph{Démonstration.} Voir \cite{Lag}, pages 455--456.\,$\square$\\

\section{Estimations finales}

\begin{theo}\label{grh}Si $\zeta_L$ vérifie l'hypothèse de Riemann généralisée, alors 
\begin{eqnarray}\left|\psi_C(x)-\frac{|C|}{|G|}x\right| &\leq & \frac{|C|}{|G|}\sqrt{x}\ln(x)\Big[\left(\frac{23}{3}+\frac{4781}{96\ln(x)}\right)\ln(d_L)\nonumber\\
& \quad + &\left(\frac{863}{31}\ln(x)+\frac{68}{3}+\frac{58681}{113\ln(x)}\right)n_L\Big]\nonumber\end{eqnarray}pour tout $x \geq 2$.\end{theo}

\emph{Remarque.} Plus précisément, on a, pour tout $x\geq 2$,\begin{eqnarray}\left|\psi_C(x)-\frac{|C|}{|G|}x\right| &\leq & \frac{|C|}{|G|}\sqrt{x}\Big[\left(\frac{23}{3}\ln(x)+\frac{29}{3}+\frac{336}{17\ln(x)}+\frac{26\ln(x)}{9\sqrt{x}}+\frac{94}{7\sqrt{x}}\right)\ln(d_L)\nonumber\\
& \quad + &\Big(\frac{863}{31}(\ln(x))^2+\frac{68}{3}\ln(x)+\frac{1198}{13}+\frac{1343}{6\ln(x)} \nonumber \\
& \quad + &\frac{77\ln(x)}{4\sqrt{x}}+\frac{3817}{30\sqrt{x}}+\frac{3}{40x^{3/4}(\ln(x))^2}\Big)n_L\Big].\nonumber\end{eqnarray}

\emph{Démonstration.} Si $\zeta_L$ vérifie l'hypothèse de Riemann généralisée, alors toutes les fonctions $\mathrm{L}(\cdot,\chi)$ associées à $\zeta_L$ par la formule $\zeta_L = \prod_\chi \mathrm{L}(\cdot,\chi)$ la vérifient également. Ainsi, pour chaque $\chi$ il n'existe pas de zéro $\rho$ non trivial tel que $|\rho|<\frac{1}{2}$, et par le lemme \ref{nchi}:

\begin{eqnarray}
\left|\sum_{\substack{\rho \\ |\gamma| < T}} \frac{x^\rho}{\rho}+\sum_{|\rho|<\frac{1}{2}} \frac{1}{\rho}\right| 
 & \leq & \sqrt{x} \left(2n_\chi(0)+\sum_{0 \leq j \leq \frac{T-1}{2}} \frac{n_\chi(2j+2)+n_\chi(-(2j+2))}{2j+1}\right) \nonumber \\
 & \leq & \frac{5}{2}\sqrt{x}\Big[\left(2+\frac{\ln(T)}{2}\right)\left(\ln(A(\chi))+\frac{1075}{134}n_E\right)\nonumber\\
 & \quad + & \left(\frac{281}{55}+\frac{\ln\left(T+4\right)\ln(T)}{4}\right)n_E\Big],\nonumber
\end{eqnarray}
car \begin{eqnarray}\int_1^{\frac{T-1}{2}} \frac{\ln(2t+5)}{2t+1}\ud t & = & \frac{1}{2}\int_1^{\frac{T-1}{2}} \left(\frac{\ln(2t+5)}{2t+1}+\frac{\ln(2t+1)}{2t+5}\right)\ud t \nonumber\\
& \quad + & \frac{1}{2}\int_1^{\frac{T-1}{2}} \left(\frac{\ln(2t+5)}{2t+1}-\frac{\ln(2t+1)}{2t+5}\right)\ud t\nonumber\\
& \leq & \frac{\ln(T)\ln(T+4)-\ln(3)\ln(7)}{4}+\frac{1}{2}\left(\frac{4}{3}\ln(2)+9\ln\left(\frac{3}{2}\right)\right).\nonumber\end{eqnarray}
En se rappelant la définition de $S(x,T)$ donnée dans le théorème \ref{formulexplicite}, ceci implique que pour tout $T \geq 2$, $|S(x,T)|$ admet la même majoration, à une multiplication par $\frac{|C|}{|G|}$ près, en remplaçant $\ln(A(\chi))$ par $\ln(d_L)$ et $n_E$ par $n_L$. 
On pose $T=\sqrt{x}+1$, par exemple (on a bien $T \geq 2$), et la formule explicite du théorème \ref{formulexplicite} prouve le résultat désiré.\,$\square$
\begin{theo}Soit $\beta_0$ le zéro réel éventuel de $\zeta_L$ vérifiant $\beta_0 \geq 1 - \frac{1}{4\ln(d_L)}$, et $\chi_0$ le caractère (réel) tel que $\mathrm{L}(\beta_0,\chi_0)=0$. Si $x \geq \exp\left(4n_L (\ln(150867d_L^{44/5}))^2\right)$, alors \[\psi_C(x) = \frac{|C|}{|G|}x - \frac{|C|}{|G|}\chi_0(g)\frac{x^{\beta_0}}{\beta_0}+\frac{|C|}{|G|}R(x),\]où $|R(x)| \leq 1505234280710x\exp\left(-\frac{7-4\sqrt{3}}{5} \sqrt{\frac{\ln(x)}{n_L}}\right)$. Le second terme peut être supprimé en l'absence du zéro exceptionnel $\beta_0$.\label{sansgrh}\end{theo}

\emph{Démonstration.} Pour alléger les calculs, posons $a =\exp\left(\frac{15}{44}\right)$, $b=3^{\frac{4}{5}}\exp\left(\frac{2156}{335}\right)$ et $c=2\frac{7-4\sqrt{3}}{5}$. Si $\rho = \beta+i\gamma$ avec $\rho \neq \beta_0$ est un zéro non trivial d'une fonction $\mathrm{L}$ tel que $|\gamma| < T$, alors la borne inconditionnelle du lemme \ref{sanszero} montre que
\[|x^\rho| = x^\beta \leq x \exp\left(-\frac{c \ln(x)}{\ln\left((ad_L)^{44/5} (b(T+3))^{n_L}\right)}\right)\]
pour $x\geq 2$ et $T \geq 2$. De plus, le lemme \ref{nchi} montre, en imitant le raisonnement de la démonstration du théorème \ref{grh}, que

\[\sum_\chi \sum_{\substack{|\rho|\geq \frac{1}{2} \\ |\gamma| \leq T}}\left|\frac{1}{\rho}\right| \leq \frac{5}{2}\left[\left(2+\frac{\ln(T)}{2}\right)\left(\ln(d_L)+\frac{1075}{134}n_L\right) + \left(\frac{281}{55}+\frac{\ln\left(T+4\right)\ln(T)}{4}\right)n_L\right].\]

Par ce même lemme, comme $\rho \neq 1-\beta_0$ implique $|\rho| \geq \frac{1}{4\ln(d_L)}$ (par le lemme \ref{presquesanszero}), on a
\begin{eqnarray}\left|\sum_\chi \sum_{\substack{\rho \neq 1-\beta_0\\ |\rho|<\frac{1}{2}}}\left(\left|\frac{x^\rho}{\rho}\right|+\left|\frac{1}{\rho}\right|\right)\right| & \leq & (\sqrt{x}+1)\sum_\chi \sum_{\substack{\rho \neq 1-\beta_0 \\ |\rho|<\frac{1}{2}}} \left|\frac{1}{\rho}\right| \nonumber \\
& \leq & 5(\sqrt{x}+1)\left[\ln(d_L)+n_L\left(\ln(3)+\frac{1075}{134}\right)\right]\ln(d_L)\nonumber\\
& \leq & \frac{8945}{9}(\sqrt{x}+1)(\ln(d_L))^2\nonumber\end{eqnarray}
par l'inégalité de Hermite-Minkowski $d_L \geq \frac{\pi}{3}\left(\frac{3\pi}{4}\right)^{{n_L}-1}$, valable quand $n_L>1$ (ce dont on déduit plus précisément $n_L \leq \frac{\ln(d_L)}{\ln\left(\frac{\pi}{3}\right)}$). 
À noter que si $n_L = 1$ (et $\ln(d_L) = 0$), alors cette inégalité est également vraie. Pour finir, 
\[\frac{x^{1-\beta_0}}{1-\beta_0}-\frac{1}{1-\beta_0} = x^\sigma\ln(x) \leq \sqrt{x}\ln(x)\]
pour un certain $\sigma \in [0,1-\beta_0]$. Tout ceci nous permet d'obtenir:

\begin{eqnarray}
\left|S(x,T) - \frac{|C|}{|G|}\chi_0(g)\frac{x^{\beta_0}}{\beta_0}\right| & = & \frac{|C|}{|G|}\left[\sum_\chi \bar{\chi}(g) \left( \sum_{\substack{|\rho|>\frac{1}{2} \\ |\gamma|<T}}\frac{x^\rho}{\rho}+\sum_{\substack{|\rho|<\frac{1}{2} \\ |\gamma|<T}}\frac{x^\rho}{\rho}-\sum_{|\rho|<\frac{1}{2}} \frac{1}{\rho}\right)-\chi_0(g)\frac{x^{\beta_0}}{\beta_0}\right]\nonumber \\
 & \leq & \frac{|C|}{|G|}\frac{5}{2}\Big[\left(2+\frac{\ln(T)}{2}\right)\left(\ln(d_L)+\frac{1075}{134}n_L\right)\nonumber\\
 & \quad + & \left(\frac{281}{55}+\frac{\ln\left(T+4\right)\ln(T)}{4}\right)n_L\Big]x\exp\left(-\frac{c \ln(x)}{\ln\left((ad_L)^{44/5} (b(T+2))^{n_L}\right)}\right) \nonumber \\
 &\quad +& \frac{|C|}{|G|}\left(\frac{8945}{9}(\sqrt{x}+1)(\ln(d_L))^2 + \frac{x^{1-\beta_0}}{1-\beta_0}-\frac{1}{1-\beta_0}\right)\nonumber \\
& \leq & \frac{|C|}{|G|}\frac{5}{2}\Big[\frac{4607}{10}+\frac{2012}{23}\ln(T)+\frac{(\ln(T+2))^2}{4\ln\left(\frac{\pi}{3}\right)}\Big]\ln(d_L)\nonumber\\
& \quad \times & x\exp\left(-\frac{c \ln(x)}{\ln\left((ad_L)^{44/5} (b(T+2))^{n_L}\right)}\right) \nonumber \\
& \quad +&\frac{|C|}{|G|}\sqrt{x}\left(\ln(x)+ \frac{8945}{9}\left(1+\frac{1}{\sqrt{x}}\right)(\ln(d_L))^2\right).\nonumber
\end{eqnarray}
En posant $T = \frac{1}{b(ad_L)^{44/5}}\exp\left(\sqrt{\frac{\ln(x)}{n_L}}\right)-3$ et en considérant les réels $x \geq 2$ tels que $\ln(x) \geq 4n_L\left(\ln(5b(ad_L)^{44/5})\right)^2$, on a bien $T \geq 2$, et
\begin{eqnarray}\left|S(x,T) - \frac{|C|}{|G|}\chi_0(g)\frac{x^{\beta_0}}{\beta_0}\right| & \leq & \frac{|C|}{|G|}\frac{5}{2}\ln(d_L)\left(\frac{4607}{10}+\frac{2012}{23}\ln(T)+\frac{(\ln(T+2))^2}{4\ln\left(\frac{\pi}{3}\right)}\right)\nonumber \\
& \quad \times & x \exp\left(-c\sqrt{\frac{\ln(x)}{n_L}}\right)\nonumber \\
& \quad + & \frac{|C|}{|G|}\left(\frac{293321}{69696}\sqrt{x}\ln(x)+\left(\frac{5}{88}\right)^2\frac{\ln(x)}{n_L}\right) \nonumber \\
& \leq & 11715979982\frac{|C|}{|G|}x\exp\left(-\frac{c}{2}\sqrt{\frac{\ln(x)}{n_L}}\right)\nonumber\\
& \quad + &\frac{|C|}{|G|}\frac{320}{19}n_Lx\exp\left(- c\sqrt{\frac{\ln(x)}{n_L}}\right),\label{majpsi1}\end{eqnarray}
d'après les inégalités (\ref{inegln1}) et (\ref{inegln2}) listées ci-après. On invoque à nouveau le théorème \ref{formulexplicite}. Ce choix pour $T$, ainsi que les inégalités suivantes,

\begin{equation}\sqrt{\frac{\ln(x)}{n_L}} \leq \exp\left(\frac{1}{a}\sqrt{\frac{\ln(x)}{n_L}}\right)+a(\ln(a)-1)\textrm{ pour tous }x \geq 2, a > 1,\label{inegln1}\end{equation}
\begin{equation}\sqrt{\frac{\ln(x)}{n_L}} \leq \exp\left(\frac{1}{4a^2n_L}\right)x\exp\left(-\frac{1}{a}\sqrt{\frac{\ln(x)}{n_L}}\right)\textrm{ pour tous }x \geq 2, a > 0,\label{inegln2}\end{equation}
\[\frac{44}{5}\ln(d_L) \leq \ln(5b(ad_L)^{44/5}) \leq \frac{1}{2}\sqrt{\frac{\ln(x)}{n_L}} \textrm{ (supposée vérifiée ci-dessus)},\]
\[n_L \leq \frac{5}{88\ln\left(\frac{\pi}{3}\right)}\sqrt{\frac{\ln(x)}{n_L}}\textrm{ (Hermite-Minkowski et l'inégalité ci-dessus)}\]
impliquent que:
\begin{eqnarray}\label{majpsi2}
\left|\psi_C(x) - \frac{|C|}{|G|}x +S(x,T)\right| & \leq & 1505243591773\frac{|C|}{|G|}x\exp\left(-\frac{1}{8}\sqrt{\frac{\ln(x)}{n_L}}\right) \nonumber\\
& \quad + & \frac{|C|}{|G|}\left[\ln(x)\left(\frac{2}{\ln(2)}+\frac{77}{4\ln\left(\frac{\pi}{3}\right)}\right)+\frac{94}{7} + \frac{1647}{5\ln\left(\frac{\pi}{3}\right)}\right]\ln(d_L)\nonumber\\
& \leq & 1505243592416\frac{|C|}{|G|}x\exp\left(-\frac{1}{8}\sqrt{\frac{\ln(x)}{n_L}}\right).
\end{eqnarray}
Finalement, le rassemblement de (\ref{majpsi1}) et (\ref{majpsi2}) donne le théorème.\,$\square$\\

Pour passer de $\psi_C$ à $\pi_C$, le raisonnement est assez classique. Posons 
\[\theta_C(x) = \sum_{\substack{\Pre \textrm{ non ramifié} \\ \mathrm{N}(\Pre) \leq x \\ \left[\frac{L/K}{\Pre}\right] = C}} \ln(\mathrm{N}_{K/\Q}(\Pre)).\]
Il y a au plus $n_K$ idéaux $\Pre^m$, avec $\Pre$ premier, dont la norme a une certaine valeur donnée, et de plus
\[
\sum_{\substack{\Pre, m \geq 2 \\ \Pre\textrm{ non ramifié}\\ \mathrm{N}_{K/\Q}(\Pre^m) \leq x}} \ln(\mathrm{N}(\Pre)) = \theta_C(x^{1/2}) + \theta_C(x^{1/3}) + \cdots + \theta_C(x^{1/[\ln(x)/\ln(2)]})
\]
car $x^{1/M} < 2$ pour $M > \ln(x)/\ln(2)$, donc l'inégalité $\theta_C(x) \leq n_K \theta_\Q(x) < 1,01624n_K x$ (tirée de \cite{Ros}, théorème 9, valable pour tout $x\geq 2$) implique que:
\begin{equation}
\psi_C(x)-\theta_C(x) = \sum_{\substack{\Pre, m \geq 2 \\ \Pre\textrm{ non ramifié}\\ \mathrm{N}(\Pre^m) \leq x}} \ln(\mathrm{N}_{K/\Q}(\Pre)) \leq \frac{22}{15}n_K\sqrt{x}\ln(x).
\end{equation}
Ceci prouve que $\theta_C$ vérifie presque la même formule asymptotique que $\psi_C$:

\begin{eqnarray}\left|\theta_C(x)-\frac{|C|}{|G|}x\right| &\leq & \frac{|C|}{|G|}\sqrt{x}\ln(x)\Big[\left(\frac{23}{3}+\frac{4781}{96\ln(x)}\right)\ln(d_L)\nonumber\\
& \quad + &\left(\frac{863}{31}\ln(x)+\frac{362}{15}+\frac{58681}{113\ln(x)}\right)n_L\Big]\nonumber\end{eqnarray}
si on suppose l'hypothèse de Riemann généralisée, et, inconditionnellement,
\[\theta_C = \frac{|C|}{|G|}x - \frac{|C|}{|G|}\chi_0(g)\frac{x^{\beta_0}}{\beta_0} + \frac{|C|}{|G|}R_0(x),\]
où $\beta_0$ est l'éventuel zéro exceptionnel de $\zeta_L$ dans la région décrite dans le lemme \ref{presquesanszero}, $\chi_0$ le caractère tel que $\mathrm{L}(\beta_0,\chi_0) = 0$, avec $|R_0(x)| \leq 1505234280719x\exp\left(-\frac{7-4\sqrt{3}}{5} \sqrt{\frac{\ln(x)}{n_L}}\right)$.

Une transformée d'Abel donne les théorèmes \ref{chebotarev} et \ref{chebotarevGRH}: on a en effet
\[\pi_C(x)-\frac{|C|}{|G|}\mathrm{Li}(x) = \frac{\theta_C(x)-\frac{|C|}{|G|}x}{\ln(x)}+\int_2^x \frac{\theta_C(t)-\frac{|C|}{|G|}t}{t(\ln(t))^2}\ud t.\]
On utilise la majoration $\int_2^x \frac{\ud t}{\sqrt{t}\ln(t)} \leq \frac{4\sqrt{x}}{\ln(x)}$, obtenue \emph{via} le changement de variable $u = \frac{-\ln(t)}{2}$ qui nous ramène à l'étude de l'exponentielle intégrale, dont on trouve toute une étude dans \cite{AS}, page 228.\\
Il faut toutefois être vigilant pour obtenir le théorème inconditionnel: comme le théorème \ref{sansgrh} ne donne une estimation valable qu'à partir de $x \geq \exp\left(4n_L (\ln(150867d_L))^2\right) =:a$, on doit faire un léger découpage:
\begin{eqnarray}\pi_C(x)-\frac{|C|}{|G|}\mathrm{Li}(x) & =& -\frac{|C|}{|G|}\chi_0(g)\frac{x^{\beta_0}}{\ln(x^{\beta_0})} + \frac{|C|}{|G|}\frac{R_0(x)}{\ln(x)}+\int_2^{\sqrt{x}} \frac{\theta_C(t)-\frac{|C|}{|G|}t}{t(\ln(t))^2}\ud t\nonumber\\
 & \quad + & \int_{\sqrt{x}}^x \frac{\theta_C(t)-\frac{|C|}{|G|}t}{t(\ln(t))^2}\ud t\Big). \nonumber\end{eqnarray}
Si, à partir de maintenant, on suppose $x \geq a^2$, alors
\begin{eqnarray}\left|\pi_C(x)-\frac{|C|}{|G|}\mathrm{Li}(x)\right| & =& \Big|-\frac{|C|}{|G|}\chi_0(g)\frac{x^{\beta_0}}{\ln(x^{\beta_0})} + \frac{|C|}{|G|}\frac{R_0(x)}{\ln(2)}+n_K\int_2^{\sqrt{x}} \frac{2,01624}{(\ln(t))^2}\ud t\nonumber\\
 & \quad + & \int_{\sqrt{x}}^x \frac{- \frac{|C|}{|G|}\chi_0(g)\frac{t^{\beta_0}}{\beta_0} + \frac{|C|}{|G|}R_0(t)}{t(\ln(t))^2}\ud t\Big|, \nonumber\end{eqnarray}
et $\int_2^{\sqrt{x}} \ud t \leq \sqrt{x} \leq \exp\left(\frac{1}{2}\left(\frac{7-4\sqrt{3}}{5}\right)^2 \right) x\exp\left(-\frac{7-4\sqrt{3}}{5} \sqrt{\frac{\ln(x)}{n_L}}\right)$, tandis que\[\int_{\sqrt{x}}^x \frac{\exp\left(-c \sqrt{\frac{\ln(t)}{n_L}}\right)}{(\ln(t))^2}\ud t \leq 4x\frac{\exp\left(-c \sqrt{\frac{\ln(x)}{2n_L}}\right)}{(\ln(x))^2}, \]d'où le théorème 1.1.\,$\square$

\end{document}